\documentclass[a4paper, leqno, intlimits]{amsart}
\usepackage{graphicx}

\usepackage{amsmath}
\usepackage{amssymb}
\usepackage{amsxtra}
\usepackage[applemac]{inputenc}
\usepackage{euscript}

\newcommand{\psh}{\mathcal{PSH}}
\newcommand{\sh}{\mathcal{SH}}
\newcommand{\C}{\mathcal{C}}

\usepackage{mathrsfs} 
\newcommand{\F}{\mathscr{F}}
\newcommand{\E}{\mathscr{E}}
\newcommand{\Pot}{\mathscr{P}}

\newcommand{\dbar}{\bar{\partial}}
\newcommand{\de}{\partial}

\newcommand{\ddc}[1]{{dd^c}\mkern .8mu #1}
\newcommand{\ma}[2][n]{(\ddc{#2})^{#1}}    

\DeclareMathOperator{\laplace}{\partial\bar{\partial}} 
\newcommand{\matwo}[1]{(\ddc{#1})^2} 

\newcommand{\bd}{\partial}

\DeclareMathOperator{\dist}{dist} 
\newcommand{\Cn}{\mathbf{C}^n}
\newcommand{\Cone}{\mathbf{C}^{1}}
\newcommand{\Ctwo}{\mathbf{C}^{2}}
\newcommand{\Real}{\mathbf{R}}
\newcommand{\Complex}{\mathbf{C}}

\def\XXint #1#2#3{{\setbox 0=\hbox{$#1{#2#3}{\int }$} 
\vcenter{\hbox{$#2#3$}}\kern -.5\wd 0}} 

\newcommand{\Proof}{\noindent\textit{Proof.}\hspace{0.6em}}
\newcommand{\QED}{\qed \vspace{\baselineskip}}

\newcommand{\weakto}{\overset{*}{\to}}
\DeclareMathOperator{\supp}{supp}

\newcommand{\ie}{i.{\kern0.11111em}e.\@}
\newcommand{\Eg}{E.{\kern0.11111em}g.\@}

\newcommand{\K}{\EuScript{K}}

\newcommand{\lp}{lp}

\newtheorem{theorem}{Theorem}[section] 
\newtheorem{proposition}[theorem]{Proposition} 
\newtheorem{corollary}[theorem]{Corollary}
\newtheorem{lemma}[theorem]{Lemma}

\theoremstyle{definition} 

\newtheorem{definition}[theorem]{Definition}  
\newtheorem{example}[theorem]{Example}

\newenvironment{remark}{\medskip\par\noindent {\it Remark\/} }{\medskip\par}
\newenvironment{question}{\medskip\par\noindent {\bf
Main Question\:}\it }{\medskip\par}

\author{Jonas Wiklund}
\address{
TFM\\
Mittuniversitetet\\
S-851 70 Sundsvall\\
Sweden\\
Fax:+42-60-148875}
\email{jonas.wiklund@miun.se}
\address{
Graduate School of Mathematics\\
Nagoya University\\
Chikusa-ku, Nagoya 464-8602\\
Japan\\
Fax: +81-52-789-5397
}
\email{wiklund@math.nagoya-u.ac.jp}

\thanks{Partially supported by JSPS grant XXXX}


\title{Pluricomplex charge at weak singularities}

\begin{document}

\begin{abstract} 
    Let \( u \) be a plurisubharmonic function, defined on a
    neighbourhood of a point \( x, \) such that the complex
    Monge-Ampère operator is well-defined on \( u.  \) Suppose also
    that \( u \) has a weak singularity, in the sense that the Lelong
    number of \( u \) at \( x \) vanish.  Is it true that the residual
    mass of the measure \( \ma{u} \) vanish at \( x \)?
    
    To our knowledge there is no known example that falsifies the
    posed question.  In this paper some partial results are obtained.
    We find that for a significant subset of plurisubharmonic
    functions with well defined Monge-Ampère mass vanishing Lelong
    number does implies vanishing residual mass of the Monge-Ampère
    measure.
\end{abstract}

\maketitle

\section{Introduction and notations}

Let us denote the plurisubharmonic functions on a domain \(
\Omega \) by \( \psh(\Omega) \) and non-positive plurisubharmonic
functions by \( \psh^-(\Omega).  \) In the same manner, subharmonic 
functions on \( \Omega \) are denoted by \( \sh(\Omega) \) and 
non-positive subharmonic functions by \( \sh^-(\Omega). \)

Perhaps the most important parameter to describe the behavior of a
plurisubharmonic function near a singularity is the so called Lelong
number.  The Lelong number of a function \( u \) at \( x \in \Cn \)
can be defined as
\begin{equation}
    \nu(u,x) = \lim_{r \to 0} \frac{1}{(2\pi)^n} \int_{\|z-x\|\leq r}
    \ddc{u} \wedge (\ddc{\log{\|z-x\|}})^{n-1},
\label{eq:lelongdef}  
\end{equation}
where we use the standard,  ``non-normalized'', differential operators
\( d:=\de + \dbar \) and \( d^c:=i(\dbar - \de). \) 
For any plurisubharmonic function \( u, \) \( \ddc{u} \) is a positive 
\( (1,1) \)-current, thus the integral make sense, and it can be shown 
that the number is bounded for any plurisubharmonic function, see
Lelong's monograph on the subject \cite{Lelong:68}.

We let \( B(r;x) \) denote the ball of radius \( r \) with center
\( x. \)  Furthermore we use the abbreviated notation \(
B(r)=B_{r}=B(r;0) \).  If we define
\[
	M(u,r,x) = \sup_{z \in B(r;x)}u(z),
\]
we have
\begin{equation} 
    \nu(u,x) = \lim_{r \to 0} \frac{M(u,r,x)}{\log{r}},
\label{eq:lelongdef2}  
\end{equation}  
\textit{confer} \cite{Avanissian:61b, Kiselman:79}.
From this identity it is is immediate  that if \( F(z) \) is a
holomorphic function \( \nu(\log\|F(z)\|,x) \) is the weight of the zero 
at \( x. \)

For any smooth function \( u \) the complex Monge-Ampère operator \(
\ma{u} \) is well-defined, but for plurisubharmonic functions in
particular, smoothness can be considerable relaxed.  By Demailly
\cite{Demailly:93} it suffices that \( u \) is plurisubharmonic and
locally bounded outside a compact subset of the domain of definition.
As one would hope, it turns out that \( \ma{\log\|z\|} = (2\pi)^n
\delta_0.  \)

For the general definition of the complex Monge-Ampère operator we
refer the reader to the papers \cite{Bedford-Taylor:76,
Bedford-Taylor:82, Demailly:93, Cegrell:energy, Cegrell:04} and the
books \cite{Klimek:book, Demailly:net}.  For some of our
results it suffices to have the Monge-Ampère operator defined on
plurisubharmonic functions bounded outside a single pole, but for the
main part of the paper the full machinery of plurisubharmonic
functions with finite pluricomplex energy, from the papers of Cegrell 
\cite{Cegrell:energy, Cegrell:04} is needed.

Since we will be using partial integration, estimates, and convergence
results in as general setting as possible, we remind you about the
definition of some of the energy-classes.  Let \( \Omega \) be a
hyperconvex domain in \( \Complex^n, \) \ie\ a domain with a continuous
plurisubharmonic function \( h \) on \( \Omega \) such that \( \{ z
\in \Omega\;|\; h(z) < c \} \) is relative compact in \( \Omega, \)
for all \( c<0.  \) The class \( \E_{0}(\Omega) \) is made up of all
negative and bounded plurisubharmonic functions \( v \) on \( \Omega,
\) such that \( \lim_{z \to \zeta} v(z) = 0, \) for all \( \zeta \in
\bd \Omega, \) and \( \int_{\Omega} (\ddc{v})^n < +\infty.  \) It is
well known that \( \ma{v} \) is well defined on bounded
plurisubharmonic functions, thus the finite mass assumption makes
sense.

A function \( u \in \psh(\Omega) \) is said to be in the class \(
\E(\Omega) \) if there, for every \( p \in \Omega \) is a
neighbourhood \( \omega \) of \( p \) and a sequence \( \{u_j\}
\subset \E_0(\Omega) \) with \( u_j \searrow u \) on \( w, \) and
subject to the total mass condition \( \sup_j \int_{\Omega} \ma{u_j}
\le +\infty.  \) If the neighbourhood \( w \) can be chosen as all of
\( \Omega \) we say that \( u \in \F(\Omega).  \) Note that for \(
\Omega \) hyperconvex subset of \( \Complex^2 \) it is known that \(
\E(\Omega) = W^{1,2}\cap\psh^-(\Omega) \) \cite{Blocki:04}.  The main
point of the class \( \E \) is that is the largest possible class where the
Monge-Ampère operator is well defined.

Often it is clear from the context which domain \( \Omega \) we use, 
or it does not matter much, in that case we often drop the reference 
to the domain from our notation. 

The comparison principle is a very strong tool in pluripotential
theory, unfortunately this principle does not hold in \( \F, \) unless
the integrability condition is further strengthen so that \(
\int_{\Omega} -u\ma{u} < + \infty \) ~\cite{Cegrell:energy}, but still
for functions in \( \F \) we have for any \( \varphi\in
\psh^-(\Omega), \) that
\begin{equation}\label{eq:comparision-in-F}  
    u \leq v \implies \int_{\Omega} \varphi \ma{u} \leq \int_{\Omega}
    \varphi \ma{v}
\end{equation}
(see \cite{AAhag:02}) which, even if it is not as strong as one would
like for the purpose of solving the Dirichlet problem for \( \ma{}, \)
will be enough for our purposes.

Given a point \( x \in \Cn, \) take an open neighborhood \( O_{x}
\) of \( x, \) and let us denote the residual mass of a measure \(
\mu \) at \( x \), i.e.\ \( \mu(O_x) - \mu(O_x\setminus \{x\}), \) with
\( \mu(\{x\}).  \) Using the Riesz decomposition formula, it is a standard
exercise in potential theory to show that for subharmonic functions in
\( \Cone \):
\begin{equation} 
    \lim_{r \to 0} \frac{M(u,r,x)}{\log{r}} = \Delta{u}(\{x\}) =
    4 \laplace{u}(\{x\}).
\label{eq:laplace}  
\end{equation}

In \( \Cn, \) \( n \geq 2, \) it is well-known that the Lelong number
is dominated by the Monge-Ampère operator in the following way:
\begin{equation}
	(2\pi \nu{(u,x)})^n \leq (\ddc{u})^n(\{x\}).
\label{eq:estimate}  
\end{equation}

Note that if \( u(z_1,z_2) = \max{\{1/k\log{|z_1|}, k^2\log{|z_2|}\}}
\), then one can show that \( \matwo{u}(0) = 4 \pi^2 k \delta_0, \) 
and since \( \nu(u,0) = 1/k, \) there can be no reverse of the inequality
in Equation~\eqref{eq:estimate} above.

As has already been pointed out by Cegrell \cite{Cegrell:04} the
inequality in Equation \eqref{eq:estimate} holds whenever \( \ma{u} \)
is well-defined, and for such functions the set \( \{z\in\Cn\;|\;
\nu(u,z) > 0 \} \) is discrete.

In this paper we try to address the following question: 
\begin{question}\label{con:main}  
    Let \( u \) be a plurisubharmonic function with well-defined
    Monge-Ampère mass \( u \in \E \), say.  Suppose the Lelong number
    of \( u \) at \( x \) is \( 0 \), is it true that \( (\ddc{u})^n
    \) does not charge the point \( x \)?
\end{question}

\begin{remark}
    First of all we note that this is purely a local problem.
    There is no difference in asking this question for functions
    in \( \E(\Omega) \) or in \( \F(\Omega), \) since if \( u \in \E(\Omega) 
    \) and \( D  \)  is an open, relatively compact subset of \( \Omega 
    \) there is a function \( u_D \in \F(\Omega) \) such that \( 
    u_D = u \) in \( D. \)
\end{remark}

There are some vague reasons to belive that, in general, \(
(\ddc{u})^n \) do not charge the point \( x \).  \Eg\ it follows
directly from a theorem in \cite{Rashkovskii:01} that if \(
u(z_{1},\ldots,z_{n}) = u(|z_{1}|,\ldots,|z_{n}|), \) \(
u^{-1}\{-\infty\} = \{0\}, \) and \( \nu(u,0) = 0 \) then \(
(\ddc{u})^n \) does not charge the origin.  A similar problem has also
been studied in for example \cite{Favre-Guedj:01}.

In this paper we prove the following two theorems

\begin{theorem}
    If \( u \geq p_{\mu}, \) where \( p_{\mu} \) is the potential of
    the pluricomplex Green function with a single pole, then vanishing
    Lelong number implies vanishing residual Monge-Ampère mass. (For
    the statement in full generality see
    Theorem~\ref{thm:potestimate}.)
\end{theorem}

\begin{theorem}\label{thm:radial2}  
    Assume that  \( u \in
    \psh \cap L^{\infty}(D^2\setminus K), \) for some \( 0 \in K \Subset \Omega.
    \) If \( u(|z|,w) = u(z,w) \) and \( \nu(u,0)=0, \) then \(
    \matwo{u}(\{0\}) = 0.\)
\end{theorem}

This article is an expanded version of a manuscript earlier published
in the authors doctoral thesis, written during a visit to Graduate
School of Mathematics, Nagoya University.  As such the author is
naturally indebted to his advisor Urban Cegrell for suggestions and
comments, and to the faculty opponent Alexander Rashkovskii for many
fruitful suggestions.  However, whatever faults the paper have is
the authors own.

\section{Some observations}

One of the main tools in analysis, and in particular in pluripotential
theory, is partial integration.  We will frequently apply partial
integration as a technique to estimate the Lelong-number.

Since the following ``Hölder like'' theorem, a rather non-trivial
application of partial integration, proved in Cegrell's seminal paper
\cite{Cegrell:04}, will be one of the main tools we state it here for
future reference.

\begin{theorem}\label{lem:holder}  
    Suppose \( u,v \in \F, \) \( h \in \E_0, \) and that \( p, q \) 
    are positive natural numbers such that \( p+q = n. \) Then
    \[ 
      \int_{\Omega} -h \, (\ddc{u})^p \wedge (\ddc{v})^q \leq
      \Big(\int_{\Omega} -h \, \ma{u} \Big)^{\frac{p}{n}}
      \Big(\int_{\Omega} -h \, \ma{v}\Big)^{\frac{q}{n}}
    \]
\end{theorem}

\Proof Cf.\ \cite{Cegrell:04}.

Let \( h_j(z) = \max{(1/j \log\|z\|, -1)}, \) then \( \{h_j\}_j \) is an
increasing sequence of continuos plurisubharmonic functions on the
unit ball of \( \Complex^n, \) tending to
\( 0 \) outside the origin.  Assuming that \( u \in
\psh{(\overline{B})}, \) we can rewrite the definition of the Lelong
number of \( u \) at the origin as
\begin{multline*}
   \lim_{r\to 0} \frac{1}{(2\pi)^n}
   \int_{B(r)} \ddc{u} \wedge \ma[n-1]{(\log\|z\|)} = \\
   = \lim_{j \to
   \infty}  \frac{1}{(2\pi)^n}
   \int_{B(1)} h_j \ddc{u} \wedge \ma[n-1]{(\log\|z\|)},
\end{multline*}
making the last integral a prime target for partial integration, or
for Theorem~\ref{lem:holder}.

The \textit{pluricomplex Green function} with \textit{pole at} \( w \)
was introduced by Klimek \cite{Klimek:85} and
Zahariuta~\cite{Zahariuta:84}.  For any connected open subset \(
\Omega \) of \( \Complex^n \) we have \( g_{\Omega}(z,w) = \sup\{u(z)\;|\;
u \in \psh^-(\Omega), \mbox{ and } \nu(u,w) \geq 1\}.  \)

Note that a continuity result by Demailly \cite{Demailly:87}, shows
that
\[ 
	g_{\Omega}(z,w) \in \C(\bar{\Omega}\times \Omega
	\setminus \{z=w\}).
\]
Thus with help of the pluricomplex Green function we can construct
increasing sequences \( \{h_j\}_j \) of continuous plurisubharmonic
functions on any hyperconvex set such that \( h_j \equiv -1 \) in a
neighbourhood of any fixed point \( x \) and such that \( h_j \nearrow
0 \) outside \( x, \) just by setting \( h_j = \max\{1/j
g_{\Omega}(z,x) , -1\}.  \) By using this sequence, it follows from
Demailly's comparison theorem of Lelong numbers \cite{Demailly:93},
that we
can express the Lelong number as
\[ 
	\nu(u,x) = \lim_{j \to \infty} \frac{1}{(2\pi)^n} \int_{\Omega}
	-h_j \, \ddc{u} \wedge \ma[n-1]{g_{\Omega}(z,x)}.
\]

Before the main body of the article we will use partial integration to
make some observations regarding the question.  First of all we note
that class of functions without concentrated mass at a point form a
convex cone, and the class is closed under taking maximum.

\begin{proposition}\label{prop:cone}  
    Let \( u,v \in \E(\Omega).  \) Assume that neither \( \ma{u} \),
    nor \( \ma{v} \) has an atom at the point \( x \in \Omega, \) then
    the same holds for \( \ma{(u+v)} \), and \(
    \ma{(\max{\{u,\varphi\}})} \) for any \( \varphi \in \E(\Omega). \)
\end{proposition}

\begin{proof}
Choose \( x \) as the origin.
Take any \( h \in \E_{0}, \) then
\[ 
	\int_{\Omega} -h\, \ma{(u+v)} = \sum_{j=0}^{n} \binom{n}{j}
	\int_{\Omega} -h\, \ma[j]{u} \wedge \ma[n-j]{v}
\]
and using the ``Hölder like'' Theorem~\ref{lem:holder}
\[ 
	\int_{\Omega} -h\, \ma{(u+v)} \leq \sum_{j=0}^{n} \binom{n}{j}
	\Big( \int_{\Omega} -h\, \ma{u} \Big)^{\frac{j}{n}} 
    \Big( \int_{\Omega} -h\, \ma{v} \Big)^{\frac{n-j}{n}}
\]
Let \( \{h_j\} \) be a sequence in \( E_0, \) such that \( h_j \equiv
-1 \) in a neighbourhood of the origin and \( h_j \nearrow 0 \) 
outside a sequence of shrinking balls \( B(r_j). \) Let \( j \to
\infty \) and we have proved the statement.

For the statement about the maximum note that \( \max(u,\varphi) \geq u \). 
By Equation~\eqref{eq:comparision-in-F} 
we have
\[
    \int -h_j\, \ma{(\max{(u,\varphi)})} \leq \int -h_j \, \ma{u},
\]
where \( h_j \) is the same sequence as above. Again, letting \( j \to 
\infty \) proves the statement.
\end{proof}

The following proposition is pretty well known, I think.  Using partial
integration---very much in the same spirit as above---we can prove the
following theorem.  I should mention that the short and elegant proof
is due to Urban Cegrell, and that many of the calculations in this
chapter uses the same idea as in this proof.

\begin{proposition}\label{prop:logestimate}  
    Let \( \Omega \) be hyperconvex and assume \( u \in \E(\Omega), \)
    \( \nu(u,0) = 0, \) and \( u(z) \geq \mbox{const} \cdot \log\|z\|.
    \) Then \( (\ddc{u})^n(\{0\}) = 0.  \)
\end{proposition}

\begin{proof}
Assume that \( u \in \F. \) Note that \( 0 \geq u \geq C \log\|z\|, \)
for some positive constant \( C. \) Take \( h \in \E_{0}, \) with \( h
\equiv -1 \) close to the origin, as above.  Then
\[ 
\begin{split}
	\int -h \, \ma{u} &= \int -u \, \ddc{h} \wedge \ma[n-1]{u} \\
   & \leq \int -C\log\|z\| \, \ddc{h} \wedge \ma[n-1]{u} \\
	& \leq \ldots \leq C^{n-1} \int -h \, \ddc{u} \wedge \ma[n-1]{(\log\|z\|)},
\end{split}   
\]
and we get in the same manner as above that \( \ma{u}(\{0\}) = 0. \)

The local statement for \( u \in \E \) follows as in the remark after
the main question.  
\end{proof}

A fundamental idea in this paper is to use estimates and approximation
\textit{from below.} Since, in general, the Monge-Ampère is only well
behaved for monotonically decreasing sequences, we need to be sure that
we can use approximation from below.
To use an approximation from below in conjunction with the Monge-Ampère
operator we use a powerful convergence theorem in \( \F. \)

\begin{theorem}\label{Thm:CX}  
    Assume that \( u_{1}^j,\ldots u_{n}^j \in \F{(\Omega)} \), \( 
    \Omega  \) hyperconvex.
    If \( u_p^j \) are monotonically increasing sequences such that \(
    u_{p}^{j} \nearrow u_{p} \) (a.e), for all \( 1\leq p\leq n, \)
    then
    \[
      \ddc{u_{1}^j}\wedge \ldots \wedge \ddc{u_{n}^j} \to
      \ddc{u_{1}}\wedge \ldots \wedge \ddc{u_{n}}, \mbox{ as } j \to 
      \infty
    \] 
    in the weak-* topology.
\end{theorem}

\begin{proof}
Since monotone convergence implies convergence in capacity, this
follows directly from Theorem 1.1 in \cite{Cegrell:01}.
\end{proof}

For decreasing sequences the convergence in weak-* sense is a
well-known property.  Bedford and Taylor showed in
\cite{Bedford-Taylor:82} that \( \ma{\cdot} \) is continuous under
\emph{decreasing} sequences in \( L^{\infty}, \) this is generalized
to decreasing sequences in \( \E \) in \cite{Cegrell:04}.

Denote the \emph{upper semicontinuous regularization} of a function \(
f \) by \( f^*.  \) Functions with one concentrated singularity, and
maximal outside this singularity is of special interest.

\begin{theorem}\label{thm:maximality}  
    Let \( \Omega \) be a hyperconvex domain containing the origin.
    If there exist a plurisubharmonic negative function \( u \in 
    \F(\Omega) \),  with \(
    \nu(u,0)=0, \) such that the Monge-Ampère operator charges the
    origin, there exist such function with \( (\ddc{u})^n
    \equiv 0 \) outside the origin.
\end{theorem}

\begin{proof}
Take \( u \in \F, \) \( r > 0, \) and define
\[ 
    S_r(z) = (\sup\{v(z)\in\psh(\Omega)\,|\, v \leq u \mbox{ on } 
    B(r); \, v \leq 0 \})^{*}.
\]

Clearly \( S_r \in \F(\Omega) \) and \( S_r \leq v \) on \( B(r). \) 
Take \( r' < r, \) then we must have that \( S_r \leq S_{r'}, \) since 
\( S_r \leq u \) on \( B(r') \) so \( S_r \) is one of the 
``competitors'' in the class of functions we take supremum of when
defining \( S_r'. \)

Let \( S = (\lim_{r \to 0} S_r )^*.  \) Then \( S \) is
plurisubharmonic and, by construction, maximal outside the origin.

Since \( S > S_r, \) we have that \( \nu(S,0) \leq \nu(S_r,0) \leq 
\nu(u,0) = 0. \)

Furthermore \( (\ddc{u})^n(\{0\}) = c, \)  \( c > 0, \) and we have
\[ 
    \int_{\Omega} (\ddc{S_r})^n \geq \int_{B_r} \ma{u} = c.
\]
Since \( S_r \) is an increasing sequence we get that \(
(\ddc{S})^n(\{0\}) \geq c, \) according to Theorem~\ref{Thm:CX}

\end{proof}

\begin{remark}
    It is important, in the construction of \( S \) above, that \(
    \ma{u} \) charges the origin.  Otherwise \( (\lim_{r \to 0} S_r)^*
    \) might be identically zero.  For example, take \( u(z) =
    -\sqrt{-\log{\|z\|}} \), then a simple calculation yields that \(
    S_r(z) = \log{\|z\|}/\log{r}, \) outside \( B(r), \) and therefore
    \( (\lim_{r\to 0}S_r)^* = 0.  \)
\end{remark}    

\section{Using estimates from below}

A radial subharmonic function whose Laplace mass does not charge the
origin can be minorized by any logarithm close to the origin. Let us
make this simple observation precise.

\begin{lemma}\label{lem:trivial}  
    Let \( D \) denote the unit disc in \( \Cone, \) and suppose \( u
    \in \sh{(D)}, u \not \equiv -\infty \) is radial (i.e.\ \(
    u(|z|)=u(z).  \) For any \( \epsilon > 0, \) let \( D_{\epsilon} =
    \{z: u(z) > \epsilon \log{|z|}\} \cup \{0\} .  \) Then \( \nu(u,0)
    = 0 = \laplace{u}(\{0\}) \) if and only if \( D_{\epsilon} \) is
    a disc of positive radius centered at the origin for all \(
    \epsilon>0.  \)
\end{lemma}

\begin{proof}
Since \( u \) is radial and \( u \not \equiv -\infty, \) it is a 
convex function in \( \log{r}, \) continuous outside of the 
(possible) pole in the origin. Assume \( \nu(u,0) = 0, \) 
Equation~\eqref{eq:lelongdef2} gives that there is a sequence \( 
(r_{j}), r_{j} \to 0 
\) such that \( u > 
\epsilon \log{r_{j}}. \)  If we change 
variables \( t \mapsto e^t, \) we get \( u(t_{j}) > \epsilon t_{j}. \) 

Suppose there is a \( t' < t_{k}, \) such that \( u(t') < \epsilon t'. 
\) Since \( u \) is convex and \( u(t_{k+1}) > \epsilon t_{k+1}, \) 
we have that \( u(t) < \epsilon t, \) for \( t < t_{k+1}, \) which
contradicts the assumption that \( \nu(u,0)=0. \)

The opposite implication follows from the continuity of the Laplace
operator.
\end{proof}

Now if we wanted to show that \( M(u,r,0)/\log{r} \) tends to zero
with \( r \) for
radial potentials which does not charge the origin, we could use the
Lemma above together with the continuity of the Laplace operator on
increasing sequences.

If we want to use this method to deal with more general plurisubharmonic
functions we need to replace the comparison with an increasing
sequence of logarithms to any sequence of plurisubharmonic functions
increasing to zero, since even for non-radial subharmonic functions on
the plane it is clearly not the case that vanishing Laplace mass
implies that we can estimate the function from below with logarithms.

To deal with more oscillating functions we mimic the sets \(
D_{\epsilon} \) in Lemma~\ref{lem:trivial} above to make the following
convenient definition.

\begin{definition}\label{def:k'}  
    Let \( U \) be an open neighbourhood of the origin.  We say that \(
    u\in \psh(U) \) is of class \( \K \) if there exists an increasing
    sequence \( \{f_j\}\subset\psh^-{(U)} \) such that \( f_j \nearrow 0
    \) (a.e.) and \( \forall j \), \( \exists r = r(j) > 0 \) such that
    \[ 
    	\{ f_j(z) \leq u( z) \} \cup \{0\} \supset B_r.
    \]
    where \( B_r \) is the ball of radius \( r, \) centered at the
    origin.
\end{definition}

It turns out that Definition~\ref{def:k'} is a handy way to describe
functions with no residual Monge-Ampère mass at the origin. The 
following theorem makes it clear.

\begin{theorem}\label{thm:class-of-k}  
    Let \( \Omega \) be an hyperconvex domain in \( \Cn \) such that \( 
    0 \in \Omega, \) and let \( u \in \F(\Omega), \) then \( u \in \K 
    \) if, and only if, \( \ma{u}(\{0\}) = 0. \)
\end{theorem}

\begin{proof}
Suppose \( u \in \K, \) then by definition there is \( f_j \in
\psh^-(\Omega) \) such that \( f_j \nearrow 0 \) and \( f_j < u \) on
balls \( B(r_j).  \)

Define a sequence of functions
\( 
	u_j := \max(u,f_j).
\)
Then \( u_j \in \F, \) and \( u_j \nearrow 0 \) (a.e.).
Clearly  \( u_j = u \) on \( B(r_j), \) thus
\[
        \int_{B(r_j)} -\varphi \, \ma{u} = \int_{B(r_j)} -\varphi \,  
        \ma{u_j}
        \leq \int_{\Omega} -\varphi \,  \ma{u_j}.
\]
Take \( \varphi \equiv -1 \) at a neighbourhood of the origin.  Since we
have, according to Theorem~\ref{Thm:CX}, \( \ma{u_j} \overset{*}{\to}
0, \) we get in particular that \( \int_{\Omega} -\varphi \ma{u_j} \to
0.  \) Thus \( \ma{u}(\{0\}) = 0.  \)

On the other hand,
suppose \( u \in \F \) such that \( \ma{u}(\{0\}) = 0. \) As in the 
proof of Theorem~\ref{thm:maximality} define a sequence
\begin{equation*}
    f_j(z) = \Big(\sup\{\varphi(z)\in \psh{(D)}\;;\; \varphi \leq 0\mbox{
    and } \varphi|_{B(1/j)} \leq u \}\Big)^*
\end{equation*}
Then \( f_j \) is an increasing sequence of plurisubharmonic functions 
such that \( f_j \nearrow 0 \) (a.e.).
\end{proof}

\begin{lemma}\label{lem:max-in-K}
   Let \( u \in \F(\Omega).  \) If \( v \in \K \) then \( \max(u,v) \in
   \K. \)
\end{lemma}

\begin{proof}
Since \( v \) is of class \( \K \) take the increasing sequence \( f_j 
\) from the definition of the class \( \K. \)  Now \( \max(u,f_j) 
\leq \max(u,v) \) on \( B({r_j}), \) and \( \max(u,f_j) \nearrow 0 \) 
pointwise (a.e.). Thus \( \max(u,f_j) \) is the required sequence for \( 
\max(u,v). \)
\end{proof}

\section{The pluricomplex potential}

Lelong~\cite{Lelong:89} has generalized the notion of the pluricomplex
Green function and defined the general \textit{multipole Green
function} on an bounded hyperconvex set \( \Omega \subset \Cn, \) with
\emph{weighted poles} \( P = \{(a_k,w_k)\}_{k=1}^{p} \) as
\begin{equation}\label{eq:multi-green-def}  
    g_{\Omega}(z,P) = \sup\{u(z)\;|\;
    u \in \psh^-(\Omega), \;
    \nu(u,w_j) \geq a_j,\: 1 \leq j \leq p\},
\end{equation} 
where the weights \( a_k \geq 0 \) and the poles \( w_k \in \Omega. \)

If \( P = \{(1,w)\} \) we write \( g_{\Omega} (z,\{(1,w)\}) = 
g_{\Omega}(z,w), \) for the pluricomplex Green function with a single 
pole at \( w  \) with weight one.

\begin{definition}\label{def:potential}  
    \cite{Carlehed:98} Let \( \mu \) be a finite, positive measure
    with support in \( \bar{\Omega}, \) where \( \Omega \) is
    a bounded domain in \( \Cn.  \) We define the
    \textit{pluricomplex potential} of \( \mu \) as
    \begin{displaymath}
	 p_{\mu}(z) = \int_{\Omega} g_{\Omega}(z,w)\,d\mu(w),
    \end{displaymath}
    and the \textit{logarithmic potential } of \( \mu \) as
    \begin{displaymath}
	\lp_{\mu}(z) = \int_{\Omega} \log\|z-w\|\,d\mu(w).
    \end{displaymath}
    
    Note that in \( \Complex^1 \) the pluricomplex potential is just
    the ordinary Green potential (or minus the ordinary
    Green potential, depending on taste).
\end{definition}

\begin{lemma}\label{lem:welldefined}  
    Let \( \Omega \) be a hyperconvex domain, and let \( \mu \) be a
    positive finite Borel measure on \( \Omega, \) with support in \(
    K \Subset \Omega.  \) Then \( p_{\mu} \) and \( \lp_{\mu} \) is in
    \( \C(\bar{\Omega} \setminus K).  \)
\end{lemma}

\begin{proof}
For the logarithmic potential this is clear.  For the pluricomplex
potential this follows from the aforementioned continuity result of
Demailly \cite{Demailly:87}.
\end{proof}

This place us in position to prove a fundamental lemma about 
pluricomplex potentials.

\begin{lemma}\label{lem:gpinF}
    Suppose \( \mu \) is a positive finite Borel measure then
    \( p_{\mu} \in \F. \)
\end{lemma}

\begin{proof}
Let \( K \Subset \Omega. \) Define \( \mu_K = \chi_K \mu \) where \( 
\chi_K \) is the characteristic function of \( K. \) Let 
\[ 
	p(z):= \int g(z,w) \, d\mu_K(w).
\]
Take \( r < \dist(K,\complement{\Omega}) \) and let \( A = 
\{z\in\Omega\;;\; \dist(z,\complement{\Omega}) < r\}. \)

Let us define
\[ 
	\alpha = \inf\{g(z,w)\;;\; z \in A,\; w\in K\}.
\]
Since, according to Lemma~\ref{lem:welldefined}, \( p \) is continuous
away from \( K, \) there is a point \( 
(z_0,w_0), \) where \(  w_0 \in K \) and \( z_0 \in A \)
 such that \( g(z_0,w_0) = \alpha. \)

For \( z \in A \)
\begin{equation}\label{eq:ineqforp}
	p(z) = \int g(z,w) \, d\mu_K(w) 
         \geq \int \alpha \, d\mu_K(w) \geq  \alpha \mu(\Omega).
\end{equation}
Since \( g(z,w_0) \in \C(\bar{\Omega}) \) we get that \( \alpha \to 0
\) as \( r \to 0, \) and thus for any \( K \Subset \Omega \) we have
\begin{equation}\label{eq:boundaryvalues} 
	\lim_{z \to \bd \Omega} p(z) = 0.
\end{equation}

Assume for simplicity that \( \mu(\Omega) = 1, \) note that since \(
\alpha < 0 \) Equation~\eqref{eq:ineqforp} gives that
\[ 
	\max(p(z),g(z,w_0)) = p(z),\text{ for } z \in A.
\]
Thus if we set \( v(z) := \max(p(z),g(z,w_0)) \) we get, by Stokes'
theorem,
\[ 
	\int \ma{p} = \int \ma{v} \leq \int \ma{g(z,w_0)} = (2\pi)^n,
\]
independent of \( K. \) This estimate and
Equation~\eqref{eq:boundaryvalues} implies that \( p_{\mu} \in \F. \)

If \( \mu(\Omega) \neq 1 \) set \( \tilde{p} = (\mu(\Omega))^{-1} p, \)
and then \( \tilde{p} \in \F, \) which implies that \( p \in \F. \)
\end{proof}


\begin{proposition}\label{prop:logpot}  
    Let \( u(z) \) be the logarithmic potential or the pluricomplex
    potential of a finite positive Borel measure with compact support on 
    a hyperconvex domain \( \Omega. \) Suppose \( \mu(\{0\}) = 0, \) then \(
    \ma{u}(\{0\}) = 0 . \)
\end{proposition}

\begin{proof}
We prove the lemma for the pluricomplex potential. The proof 
is similar for the logarithmic potential.
 
Let \( \mu_j = \chi_j \mu, \) where \( \chi_j \) is the characteristic
function for \( B(1/j), \) and define an increasing sequence in \(
\F(\Omega) \) by the formula
\[ 
	f_j(z) = 2 \int_{\Omega} g(z,w)\,d\mu_j + (1/j) \log{\|z\|}.
\] 
Since \( \mu(\{0\}) = 0, \) we have \( f_j \nearrow 0 \) pointwise
(q.e.). Theorem~\ref{thm:potinF} guarantees that \( \{f_j\}_j \subset 
\F. \) 

Now, consider the set \( \{ u - f_j > 0\}. \)
\[
\begin{split}
(u - f_j)(z) &= \int_{\Omega} g(z,w)\,d\mu(w) - \int_{\Omega}
2g(z,w)\,d\mu_j(w)
                - \frac{1}{j}\log{\|z\|} \\
             &= \int_{\complement{B({1/j})}} g(z,w)\,d\mu(w)
                - \int_{B(1/j)} g(z,w)\,d\mu_j(w)
		- \frac{1}{j}\log{\|z\|}\\
             &= h(z) - \int_{B(1/j)} g(z,w)\,d\mu_j(w) 
	        - \frac{1}{j}\log{\|z\|},
\end{split}
\]
where \( h(z) := \int_{\complement{B({1/j})}} g(z,w)\,d\mu(w). \)

According to Lemma~\ref{lem:welldefined} \( h(z) \) is continuous on \(
B(1/j).  \) In particular \( h(z) > -M \) for some  \(  M > 0 \) on \(
{B(1/(j+1))},  \) and in addition we have 
\(
	-\int_{B(1/j)} g(z,w)\,d\mu_j(w) \geq 0,
\) on the same set.
Thus
\[ 
	(u - f_j)(z) > -M -\int_{B(1/j)} g(z,w)\,d\mu_j(w) - 1/j \log{\|z\|} > 
    -M  - 1/j \log{\|z\|} > 
    0,
\]
as long as \( \|z\| < e^{-jM}. \) 

Hence, for every \( j \), there is a \( r, \) only depending on \( j,
\) such that if \( z \in B_r \) then \( u(z)-f_j(z) > 0, \) thus \( u
\in \K \) (see Definition~\ref{def:k'}) and according to
Theorem~\ref{thm:class-of-k} we have \( \ma{u}(\{0\}) = 0.  \)
\end{proof}


\begin{corollary}\label{cor:lelong-number-of-pot}
    Let \( \mu \) be a positive finite Borel measure on \( \Omega. \) 
    Take \( x \in \Omega \) then \( \nu(p_{\mu},x) = \mu(\{x\}). \)
\end{corollary}

\begin{proof}
We may decompose Green potential of \( \mu \) as \( p_{\mu}(z) =
\int_{\Omega} g(z,w)\, d\mu(w) = \mu(\{x\}) g(z,x) + p_{\mu'}(z), \)
where \( \mu' \) does not have an atom at \( x.  \) According to
Proposition~\ref{prop:logpot} the Monge-Ampère measure of \( p_{\mu}
\) does not charge \( x, \) and since the Lelong number is dominated
by the Monge-Ampère charge we get \( \nu(p_{\mu'},x) = 0.  \)

On the other hand, by definition of the Green function, 
\( \nu(\mu(\{x\})g(z,x),x) = \mu(\{x\}).  \)
\end{proof}

\begin{definition}
    Let \( \Omega \) be a hyperconvex domain, let \( g_{\Omega}(z,w) \)
    be the pluricomplex Green function with pole at \( w.  \) 
    We define the
    class \( \Pot \) by saying that \( \varphi \in \Pot(\Omega) \) if
    \( \varphi \in \psh^-(\Omega) \) and
    there is a finite positive Borel measure on \( \Omega \) such that
    \[ 
      \varphi(z) \geq \int_{\Omega} g_{\Omega}(z,\zeta) \,
      d\mu(\zeta).
    \]
\end{definition}

\begin{theorem}\label{thm:potinF}
    Suppose \( \Omega \)  is a hyperconvex domain in \( \Complex^n, \) 
    for dimension \( n \geq 2,  \) then
    \( \Pot(\Omega) \subsetneq \F(\Omega). \)
\end{theorem}

\begin{proof} 
According to  Lemma~\ref{lem:gpinF} every Green potential is in \( \F, \)
and since any function \( u \in \Pot \) vanishes on the boundary and
is minorized by a function in \( \F \) the inclusion is clear. 

To see that the set \( \F(\Omega) \setminus \Pot(\Omega) \) is not
empty if \( \Omega \subset \Cn \) for \( n \geq 2, \) consider the
sequence of functions \( u_N(z) = g_{\Omega} (z,P_N), \) where the
weighted poles for the pluricomplex Green function is \( P_N = \{
(k^{-p/n}, w_k)\}_{k=1}^{k=N},\) for \( 1<p<2 \) fixed.

We choose the poles \( (w_k) \) such that \( \bigcup\{w_k\} \Subset
\Omega.  \) By construction of the Green function \( u_N \in
\F(\Omega) \) and
\[ 
   \int_{\Omega} \ma{u_N} = (2\pi)^n \sum_{k=1}^{N} (k^{-p/n})^n < + 
   \infty,
\]
since \( p>1.  \) Thus if we define \( u := \lim_{N\to \infty} u_N \)
we have \( u \in \F(\Omega).  \)
    
However \( u \not\in \Pot, \) because if \( u_N \geq \int g \, d\mu \)
we have 
\[ 
	k^{-p/n} = \nu(u_N, \omega_k) \leq \nu(p_{\mu}, \omega_k) =
   \mu(\{\omega_k\}),
\]
where the last equality follows from
Corollary~\ref{cor:lelong-number-of-pot}. Then
\[ 
      \mu(\Omega) \geq \sum_{k=1}^N k^{-p/n} \to +\infty, \mbox{ as } N 
      \to \infty,\mbox{ for } n \geq 2, 
\]
hence \( u \not\in \Pot. \)

Note that since functions in \( \F \) has harmonic majorant \( 0, \)
it follows from the Riesz representation theorem that \( \F = \Pot \)
in \( \Cone.  \)
\end{proof}

It is possible to reshape Proposition~\ref{prop:logpot} a bit for the purpose
of generalizing Proposition~\ref{prop:logestimate}.  We begin with a
proposition that generalize Proposition~\ref{prop:cone} and is
interesting in its own right.

\begin{proposition}\label{prop:suminF}  
    Let \( \Omega \) be a hyperconvex domain and take \( x \in \Omega. \)
    If \( f_j \in \E(\Omega), \) \( \ma{f_j}(\{x\}) = 0, \) and \(
    \sum_1^{\infty}{f_j} \in \E(\Omega), \) then \(
    \ma{(\sum_1^{\infty}{f_j})}(\{x\}) = 0.  \)
\end{proposition}

\begin{proof}
Since this is a local statement, we might as well suppose that
\( \Omega = B \) and take \( x \) as the origin.  Furthermore we might
as well assume---after modification outside a neighbourhood of the
origin and then translation---that \( f_j \in \F \) and \(
\sum_1^{\infty}{f_j} \in \E(B).  \) Take \( h \in \E_0.  \) By partial
integration:
\[ 
	 \int_B
	- (\sum_{j=1}^{\infty} f_j) \, \ddc{h} \wedge
	\big(\ddc(\sum_{j=1}^{\infty} f_j)\big)^{n-1}
   = \int_B -h \, \big(\ddc(\sum_{j=1}^{\infty} f_j)\big)^n < \infty
\]
Now, 
\( 
	\mu = \ddc{h} \wedge \big(\ddc(\sum_{j=1}^{\infty} f_j)\big)^{n-1}
\)
is a positive finite measure on \( \Omega, \) vanishing on pluripolar
sets.  Hence
\( 
	\int \sum_{j=1}^{\infty} f_j > -\infty
\)
and so
\begin{equation}
    \lim_{N \to \infty}  \sum_{j=N}^{\infty} f_j = 0 \mbox{ a.e.\ }(\mu).
    \label{eq:finitesum}  
\end{equation} 

Let us define \( h_k(z) = \max{(-1,k^{-1}\log{\|z\|})}, \) then \(
\{h_k\}_{k=1}^{\infty} \) is an increasing sequence of negative
functions, thus \( -h_1\geq \ldots \geq -h_k \geq \ldots \)
By Theorem~\ref{lem:holder} we have
\begin{eqnarray*}
    \lefteqn{ \int -h_k \,  \big(\ddc(\sum_{j=1}^{\infty} f_j)\big)^n
    } \\
	& = & \int -h_k \, \ddc{(\sum_{j=1}^{N-1} f_j)} \wedge
	\big(\ddc(\sum_{j=1}^{\infty} f_j)\big)^{n-1} + \\
    &  & \quad +
    \int -h_k \, \ddc{(\sum_{j=N}^{\infty} f_j)} \wedge
	\big(\ddc(\sum_{j=1}^{\infty} f_j)\big)^{n-1}\\
   & \leq & \int -h_k \, \ddc{(\sum_{j=1}^{N-1} f_j)} \wedge
   \big(\ddc(\sum_{j=1}^{\infty} f_j)\big)^{n-1} + \\
    &  & \quad +
    \int -h_1 \, \ddc{(\sum_{j=N}^{\infty} f_j)} \wedge
   \big(\ddc(\sum_{j=1}^{\infty} f_j)\big)^{n-1}\\
   & \leq & \Big[ \sum_{j=1}^{N-1} \big( \int -h \, \ma{f_j} \big)^{1/n}
	\Big] \Big[ \int -h \,  \big( \ddc{(\sum_{j=1}^{\infty} f_j)} \big)
	\Big]^{\frac{n-1}{n}} + \\
	 & & \quad + \int - \big( \sum_{j=N}^{\infty} f_j \big)\,
	 \ddc{h} \wedge \big( \ddc(\sum_{j=1}^{\infty} f_j)
	 \big)^{n-1}.  \\
\end{eqnarray*}

By Equation~\eqref{eq:finitesum} the last term in the estimate above
can be made arbitrarily small if \( N \) is chosen large enough.  By
assumption, \( \int_{\{0\}} \ma{f_j} = 0,\) so if we then choose \( k
\) big enough \( \int -h_k \, \ma{f_j} \) can also be made as small as
we like, and we have
\[ 
	0 = \lim_{k\to\infty} \int -h_k \, \big(\ddc(\sum_{j=1}^{\infty}
	f_j)\big)^n = \big(\ddc(\sum_{j=1}^{\infty} f_j)\big)^n(\{0\}),
\]
which proves the theorem. 
\end{proof}

\begin{remark}
    Note that the condition \( \ma{(f_j)}(\{0\}) = 0, \forall j \) does
    not suffice to guarantee that \( \sum{f_j} \in \F. \) Take for
    example \( f_j = g(z,b_j), \) where \( \{b_j\}_1^{\infty} \) is a
    sequence such that \( \{b_j\} \Subset \Omega \) and \( b_j \neq 0.
    \) Then \( \ma{f_j} (\{0\}) = 0, \) but because \( \int
    \ma{(\sum_{j=1}^N{f_j})}= (2\pi)^n N, \) we have \( \sum f_j \not
    \in \F. \)
\end{remark}   


\begin{proposition}\label{prop:nullpot}  
    Suppose \( \mu \) is a finite positive Borel measure on a
    hyperconvex domain \( \Omega.  \) Take \( x \in \Omega, \) if 
    \[
   	u(z) = \int_{\Omega} g(z,w) \, d\mu(w), 
    \] 
    and \( \nu(u,x) = 0 \) then \( \ma{u}(\{x\}) = 0.  \)
\end{proposition}

\begin{proof}
By Corollary~\ref{cor:lelong-number-of-pot}, \( \nu(u,x) = 0 \) implies
that \( \mu(\{x\} = 0, \) so this is a direct consequence of
Proposition~\ref{prop:logpot}.

Let us demonstrate how it also follows from
Proposition~\ref{prop:suminF}.  After scaling we may as well assume
that \( B(1) \subset \Omega, \) and for convenience we take \(
x \) as the origin.

If \( \mu(\{0\}) = 0, \) then
\[ 
	u(z) = \sum_{j =1}^{\infty} \int g(z,w)\, d\mu_j(w),
\]
where \( \mu_j = \chi_j d\mu, \) for \( j=0,1,\ldots, \) \( \chi_0 =
\Omega \setminus B_1, \) and \( \chi_j \) is the characteristic
function for \( \overline{B(1/j)} \setminus B(1/(j+1)), \) for \(
j=1,2,\ldots \) Then
\[ 
	\Big(\ddc{(\int g(z,w)\, d\mu_j(w))}\Big)^n(\{0\}) = 0,
\]
and according to Proposition~\ref{prop:suminF} \( \ma{u}(\{0\}). \) 
\end{proof} 


Recall that Proposition~\ref{prop:logestimate} stated that for
plurisubharmonic
functions that can be minorized by a logarithm vanishing Lelong number
implies vanishing residual Monge-Ampère mass. We are now in position
to generalize that to plurisubharmonic functions minorized by a more
general class of  functions.

\begin{theorem}\label{thm:potestimate}  
     Let \( \Omega \) be a hyperconvex domain and suppose \( u \in
     \Pot(\Omega).  \) Given a point \( x \in \Omega \) and
     \( \varphi \in \E(\Omega) \) such that \( \ma{\varphi}(\{x\}) = 0, \) 
     let \( v \geq u + \varphi, \)
     then \( \nu(v,x) =
     0 \) implies that \( \ma{v}(\{x\}) = 0.  \)
\end{theorem}

\begin{proof}
This is entirely a local statement so we might as well suppose 
that \( x = 0 \) and \( \Omega = B(0,1). \)

First we assume \( u \in \F, \) \( \nu(u,0) = 0 \) and that \( u \geq
p_{\mu} \) for some positive measure \( \mu.  \) There is a number \(
a \geq 0 \) such that \( p_{\mu}(z) = a \log{\|z\|} + s(z), \) where
\( s(z) = \int_{B\setminus\{0\}} g(z,w) \, d\mu, \) hence for any \(
h_k \in \E_0: \)
\begin{eqnarray}
    0 & \geq & \int h_k \ma{u} = \int u \, \ddc{h_k} \wedge
    (\ddc{u})^{n-1}\notag \\
    & \geq & \int u \, \ddc{h_k} \wedge (\ddc{f})^{n-1} \geq \int u \,
    \ddc{h_k} \wedge (\ddc{(a \log{\|z\|} + s(z))})^{n-1} \notag \\
    & = & \int u \, \ddc{h_k} \wedge (\ddc{a \log{\|z\|}})^{n-1} + \notag \\
    & & \; + \sum_{j=0}^{n-2} \binom{n-1}{j} \int u \, \ddc{h_k} \wedge
    (\ddc{a \log{\|z\|}})^{j} \wedge (\ddc{s(z)})^{n-1-j}.
    \label{eqarray:line4}  
\end{eqnarray} 

Let \( h_k = \max{(\log{\|z\|}/k,-1}).  \) The first term in
Equation~\eqref{eqarray:line4} tends to zero as \( k \) tends to
infinity as in Proposition~\ref{prop:logestimate}.  For the other terms we
have, according to Theorem~\ref{lem:holder},
\begin{eqnarray}
    \lefteqn{\int - u \, \ddc{h_k} \wedge (\ddc{a 
          \log{\|z\|}})^{j} \wedge (\ddc{s(z)})^{n-1-j}}\notag \\
      & \leq & \int - h_k \ddc{u} \wedge (\ddc{a 
          \log{\|z\|}})^{j} \wedge (\ddc{s(z)})^{n-1-j}\notag \\
      & \leq & 
      ( 2\pi a )^{j}
      \Big[ \int \ma{u} \Big]^{\frac{1}{n}}
      \Big[ \int -h_k \ma{s} \Big]^{\frac{n-j-1}{n}}, \label{eq:above}
\end{eqnarray}
which---according to Proposition~\ref{prop:nullpot}---tends to zero as
\( k \) tends to infinity.

Now for the general case when \( v \geq \int p_{\mu} + \varphi, \)
for some \( \varphi \in \E(\Omega) \) with \( \ma{\varphi}(\{x\}) = 0,
\) where we proceed much in the same manner.  Again there is a positive
number \( a \) such that \( p_{\mu} = a \log{\|z\|} + s(z), \)
where \( s(z) = p_{\mu'} \) with \( \mu'(\{0\}) = 0.  \) Set \(
\tilde{\varphi} = s + \varphi, \) then \( v(z) \geq a \log{\|z\|} +
\tilde{\varphi}(z), \) where, according to Proposition~\ref{prop:cone},
\( \tilde{\varphi} \in \E(\Omega) \) has the property \(
\ma{\tilde{\varphi}}(\{0\}) = 0.  \)

As in  the remark  following the main question we might as 
well assume that \( \tilde{\varphi} \in \F. \) Take \( h \in \E_0 \) then
\begin{eqnarray*}
    \lefteqn{ \int -h \,  \ma{v} }\\
     & \leq & \int -h \,  \ddc{v} \wedge 
                               \ma[n-1]{(a\log{\|z\|} + \tilde{\varphi})} \\
     & = & \int -h \,  \ddc{v} \wedge \ma[n-1]{a \log{\|z\|}} + \\
     &  & \quad + \sum_{j=1}^{n-1} \binom{n-1}{j} 
          \int -h \,  \ddc{v} \wedge \ma[n-1-j]{a \log{\|z\|}} 
                                     \wedge \ma[j]{\tilde{\varphi}}.
\end{eqnarray*}
Choosing \( h = \max(1/k \log{\|z\|},-1) \), the first term in the sum above is
arbitrarily close to zero as in Proposition~\ref{prop:logestimate}.  The
remaining terms can be taken cared of in the same way as in
Equation~\eqref{eq:above} above.  
\end{proof}

\section{Examples of functions with no Monge-Ampère mass at the
poles} \label{sec:examples} 

A function \( u:\Cn\to\Complex^p \) is said to be \textit{poly-radial}
if it is a radial function in every variable separately \(
u(|z_1|,\ldots,|z_n|) = u(z_1,\ldots,z_n).  \)

If \( u \) is a poly-radial plurisubharmonic function in a neighbourhood
of the origin it is a radial subharmonic function along any complex line
through the origin.  Since radial subharmonic functions are continuos
outside the origin, or identically \( -\infty, \) it follows that if \(
u \) is bounded below away from the origin it has its only possible
pole at the origin.

We will need a couple of lemmas about the Lelong number along 
slices. Given a function \( u:\Complex^n \supset \Omega \to \Real 
\cup \{-\infty\} \) we define a \textit{slice} of \( u \) through \( 0 
\) and \( 
p \in \Complex^n, \) \( u_{p}, \) as: \( u_{p}(\zeta):= u(\zeta p), 
\zeta \in \Complex, \) wherever this expression make sense.

\begin{lemma}\label{lem:0}  
    \( \nu(u_{p}, 0) \geq \nu(u,0) \)
\end{lemma}

\begin{proof}
Since \( \log{r} < 0, \) if \( r < 1 \) we have:
\[ 
	\lim_{r \to 0}\frac{ \sup_{\|\zeta p\| \leq r} u(z)}{\log{r}} \leq 
	\lim_{r \to 0} \frac{\sup_{\|\zeta\| \leq r} u(\zeta p)}{\log{r}}.
\]
\end{proof}

Using this lemma one can prove that the reverse inequality holds almost 
everywhere.

\begin{lemma}\label{lem:1}  
    Assume \( u \in \psh(B), \) take \( q \in B \) fixed, then \( 
    \nu(u_{q},0) = \nu(u,0) \) for all \( q \in B \setminus A ,\) 
    where \( A \) is a pluripolar set.
\end{lemma}    

\Proof This is well known, Cf.\ \cite{Cegrell-Thorbiornson}, or for a 
more elegant proof, the survey \cite{Rashkovskii:overview}.

For functions radial in at least one variable Lemma~\ref{lem:1} above
can be considerable strengthen.  Namely, if the Lelong number at the
origin vanish, it vanish on all slices through the origin, except
perhaps along the coordinate axes.

\begin{lemma}\label{lem:2}  
    Assume that \( u \in \psh{(B)}, \) where \( B \) is the unit-ball
    in \( \Ctwo, \) and that \( u(|z|,w)=u(z,w).  \) Suppose that \(
    \nu(u,0) = 0, \) then for all \( y=(y_{1}, y_{2}) \in B, \) such
    that \( y_{1}y_{2} \neq 0, \) \( \nu(u_{y},0) = 0.  \)
\end{lemma}

\begin{proof}
Take \( p=(z,w), q=(z',w') \in \Complex^2, \) with \( |z|<|z'|, 
\) and \( |w|=|w'|=R. \) Since \( u(r,w) \) is a increasing function in 
the radial variable we have that
\[ 
    \sup_{|w|=R} u(|z|,w) \leq \sup_{|w'|=R} u(|z'|,w').
\]
That is, we have \( \nu(u_{p},0) \geq c \nu(u_q,0), \) for some 
constant. Take a point \( y=(y_{1},y_{2}) \) such that neither \( 
y_{1}, \) nor \( y_{2} \) is zero. According to Lemma~\ref{lem:1} 
there is a point \( y'=(y'_{1},y_{2}) \) with \( |y'_{1}|<|y_{1}| \) 
such that \( \nu(u_{y},0) = 0 \), but then \( 0 \geq \nu(u_{y},0). \)
\end{proof}

Now let us demonstrate how we can use the convergence on increasing
sequences (Theorem~\ref{Thm:CX}).  Using Lemma~\ref{lem:trivial} and
Theorem~\ref{Thm:CX} we can deal with poly-radial functions.  
It can be illustrative to see that how this
follows directly from the convergence Theorem~\ref{Thm:CX}.

\begin{proposition}\label{prop:poly-radial}  
    Assume \( u \in \psh(D^n), \) is poly-radial, and \(
    u^{-1}\{-\infty\} = \{0\}.  \) Then \( \nu(u,0) = 0 \) if, and
    only if, \( (\ddc{u})^n(\{0\}) = 0.  \)
\end{proposition}

\begin{proof}
Let \( \pmb{e_j} \) be the unit vector in \( \Cn \) having \(
1+i \cdot 0 \) as its \( j \):th coordinate.  Note that \(
u_{\pmb{e_j}} \) is a subharmonic function on the unit disc in \(
\Complex, \) and since \( u^{-1}\{-\infty\} = \{0\} \) it follows that
\( u_{\pmb{e_j}} \not \equiv -\infty \) and thus \(
\nu(u_{\pmb{e_j}},0) = \nu_j < +\infty.  \)

Take \( \epsilon >0  \) and define
\[ 
	v_{\epsilon}(z) := \max((\nu_1+\epsilon)
	\log|z_1|,\ldots,(\nu_n+\epsilon) \log|z_n|).
\]
Then \( v_{\epsilon} < u \) on a polydisc centered at the origin (as 
in the proof of Theorem 4.2~\cite{Wiklund:04}). Applying 
Proposition~\ref{prop:logestimate} we see that \( \nu(u,0) = 0 \) implies 
that \( \ma{u}(\{0\}) = 0. \)
\end{proof}

Proposition~\ref{prop:poly-radial} also follows more or less directly
from Demailly's comparison principle for generalized Lelong numbers
\cite{Demailly:93},  or directly from a general theorem in
\cite{Rashkovskii:01}.

The goal is to generalize the idea of the  proof of 
Proposition~\ref{prop:poly-radial} to hold for a more general class
of plurisubharmonic functions.  But the main idea is still to compare
the function with a smaller function along a line and then to amplify
the comparison to a larger domain.


To highlight the methods that are used through this section we start
off with a class of functions for which the the answer to the main
question is ~\ref{con:main}
is rather easily seen to be affirmative.

\begin{proposition}\label{prop:2}  
    Assume that  \( u \in
    \F(D^2).  \) Let \( \pmb{e_2} = (0,1+0\cdot i).  \) If \(
    \nu(u_{\pmb{e_2}},0) = 0, \) i.e.\ the Lelong number of \( u \)
    along the \( z_2 \)-axis is zero, and for any point \( z_2 \) \(
    u(0,z_2) \leq u(z_1,z_2), \) for all points \( z_1\in D, \) then \(
    \matwo{u} \) does not charge the origin.
\end{proposition}

\begin{proof}
Since the Lelong number of \( u \) along the \( z_2 \)-axis
vanish we have, after changing \( u \) near the boundary if necessary,
that \( u(0,z_2) \) is of class \( \K(D), \) thus we can apply
Lemma~\ref{lem:max-in-K} to see that \(
\max(u(z_1,z_2),u(0,z_2)) \in \K. \) Since we assumed that \( u(0,z_2)
\leq u(z_1,z_2) \) we have that \( \max(u(z_1,z_2),u(0,z_2)) =
u(z_1,z_2), \) thus \( u \in \K. \) According to
Theorem~\ref{thm:class-of-k} \( \ma{u}(\{0\}) = 0.  \) 
\end{proof}

Note that according to Lemma~\ref{lem:0} we have that functions 
satisfying the conditions in Proposition~\ref{prop:2} must have 
Lelong number zero at the origin.

For functions that are radial in at least one of the variables more 
could be said. We start off with an immediate corollary to 
Proposition~\ref{prop:2}

\begin{corollary}\label{cor:2}  
    Assume that  \( u \in 
    \F(D^2) 
    \) Take \( \pmb{e_2} = (0,1+0\cdot i). \)
    If \( \nu(u_{\pmb{e_2}},0) = 0, \) 
    and \( u(|z_{1}|,z_2) = u(z_1,z_2) \) then
    \( \matwo{u} \) does not charge the origin.
\end{corollary}

\begin{proof}
Using the maximum principle in the first variable we have \(
u(0,z_{2}) \leq u(z_{1}, z_{2}) \) and we may apply
Proposition~\ref{prop:2}.  
\end{proof}


\begin{theorem}\label{thm:maxinF}
    Let \( \Omega_1 \) and \( \Omega_2 \) be hyperconvex domains in \(
    \Complex^{n_1} \) and \( \Complex^{n_2} \) respectively. Suppose \( u_1 \in
    \F(\Omega_1) \) and \( u_2 \in \F(\Omega_2), \) then \( 
    \max(u_1,u_2) 
    \in \F(\Omega_1 \times \Omega_2) \) and
    \begin{equation}\label{eq:total-energy} 
    	\int_{\Omega_1 \times \Omega_2} \ma[n_1+n_2]{(\max(u_1,u_2))} =
        \int_{\Omega_1} \ma[n_1]{u_1} \: \int_{\Omega_2} \ma[n_2]{u_2}.
    \end{equation}
\end{theorem}

\begin{proof} 
It is enough to prove a special case of the statement, namely if
\( u_i \in \E_0(\Omega_i)\cap\C(\bar{\Omega}_i) \) 
with
\( \supp\{\ma[n_i]{u_i}\} \Subset \Omega_i, \) 
for \( i = 1,2, \)
then 
\( \max(u_1,u_2) \in \E_0(\Omega_1 \times \Omega_2) \) and 
Equation~\eqref{eq:total-energy} holds, because according to the main 
approximation Theorem 2.1 in \cite{Cegrell:04} there are always
sequences
\( u_i^j \in \E_0(\Omega_i) \) with \( \supp\{\ma[n_i]{u_i^j} \} \Subset 
\Omega_i, \) for \( i=1,2, \) 
such that
\( u_1^j \searrow u_1 \) and \( u_2^j \searrow u_2, \) as \( j \to 
\infty. \)
Hence 
\[ \lim_{j\to \infty} \max(u_1^j,u_2^j) =  \max(u_1,u_2) \in 
\F(\Omega_1 \times \Omega_2)  \] by definition.

Assume therefore that 
\begin{equation}\label{eq:compact support}
	\supp\{\ma[n_i]{u_i}\} \Subset \{ z \in \Omega_i\;;\; u_i(z) < -\delta 
	\}, \; i=1,2,
\end{equation}
for some \( \delta > 0. \)
Then
\[ 
	u_{\delta} =  \max(u_1,u_2, - \delta) =
    \max{\big\{ \max(u_1, - \delta) + \delta , \max(u_2, - \delta) + 
    \delta \big\}} - \delta
\]
so 
\[ 
	u_{\delta} + \delta = \max{\big\{ \max(u_1, - \delta) + \delta , \max(u_2, -
	\delta) + \delta) \big\}}.
\] 
Set \( \phi_i = \max(u_i,-\delta)+\delta, \) for \( i=1,2, \) then
\[ 
	\int_{\{\phi_i > 0\} } \ma[n_i]{\phi_i} = 0
\]
by the assumption of the support of the Monge-Ampère masses in
Equation~\eqref{eq:compact support}. Thus we can 
apply a theorem by B{\l}ocki (Theorem 7, \cite{Blocki:00}) to get
\[ 
	\ma[n_1+n_2]{(u_{\delta}+\delta)} = 
    \ma[n_1]{( \max(u_1, - \delta) )} \wedge
    \ma[n_2]{( \max(u_2, - \delta) )}.
\]
To sum up, we have
\begin{multline*} 
	\int_{\Omega_1 \times \Omega_2} \ma[n_1+n_2]{u_{\delta}} =
    \int_{\Omega_1} \ma[n_1]{ ( \max(u_1, - \delta) )} \: 
    \int_{\Omega_2} \ma[n_2]{( \max(u_2, - \delta) )}\\
     =
    \int_{\Omega_1} \ma[n_1]{ u_1 } \: 
    \int_{\Omega_2} \ma[n_2]{ u_2 }.
\end{multline*}
Where the last equality follows from Stokes' theorem since \( 
\max(u_i, -\delta) = u_i \) close to the boundaries.

Set \( u := \max(u_1,u_2). \)
Since \( u = \max(u_1,u_2) = \max(u_1,u_2,-\delta) = u_{\delta} \) 
outside a compact subset of \( \Omega_1 \times \Omega_2 \) we also get
\[ 
	\int_{\Omega_1 \times \Omega_2} \ma[n_1+n_2]{u_{\delta}} =
    \int_{\Omega_1 \times \Omega_2} \ma[n_1+n_2]{u} =
    \int_{\Omega_1} \ma[n_1]{ u_1 } \: 
    \int_{\Omega_2} \ma[n_2]{ u_2 }.
\]
Which proves the  required special case.
\end{proof}

\begin{example}
    Fix a number \( 0 < r < 1, \) and let \( u(z,w) =
    \max{\{\log{|z|}, -(-\log|w|)^a \}} + \log{r}, \) for \( 0<a<1.  \)
    A direct computation gives that
    \[ 
    	\int_{D_{\rho}} - \ddc{(-\log{|w|})^a} = \int_0^{\rho}
      \frac{a(1-a)(-\log{r})^{a-2}}{r} dr =
      \frac{a}{(-\log{\rho})^{1-a}} .
    \]
    Let \( \rho \) be the solution of \(
    -(-\log{\rho})^a=\log{r}. \) By Theorem~\ref{thm:maxinF}
    we get that \( u \in \F( D_r \times D_{\rho}) \) with total mass
    \[ 
    	\int_{D_r \times D_{\rho}} \ma[2]{u} = \int_{D_r} \ddc{\log|z|} 
    	\int_{D_{\rho}} \ddc{-(-\log{|z|})^a} = \frac{a}{(-\log{\rho})^{1-a}}.
    \]
    Thus \( u \in \E \setminus \F(D_1 \times D_1), \) but we may change
    \( u \) outside a neighbourhood of the origin so that \( u \)
    retains all its properties at the origin but still remain in \( \F. \)
\end{example}

\begin{example}
    More general, let \( u(z,w) = \max{\{\log{|z|}, -({-v(w)})^a \}},
    \) for \( 0<a<1, \) where \( v \in \sh^-{(D)}, \) \( v \not \equiv
    -\infty .\) According to Theorem~\ref{thm:maxinF} \( v \in \F \)
    after we have modified the function close to the boundary.
    Since  \( \nu(-({-v(w)})^a,0) = 0 \) we have \(
    \nu(u_{\pmb{e_2}},0) = 0 \) and thus, according to
    Corollary~\ref{cor:2}, \( \matwo{u} (\{0\}) = 0.  \)
\end{example}


\begin{example}
    Take a sequence of positive numbers \( \{a_{j}\}_{j=0}^{\infty}
    \subset \Real, \) such that \( \sum_{j=0}^{\infty} a_{j} <\infty,
    \) and a sequence \( \{b_{j}\}_{j=0}^{\infty} \subset \Complex, \)
    with \( b_j \neq 0 \) and \( \lim_{j\to \infty} b_{j} = 0.  \) Let
    \[ 
    	u(z,w)=\max\{\sum_{j=0}^{\infty} a_{j}\log{|z-b_{j}|}, \log|w|
    	\},
    \] 
    then according to Theorem~\ref{thm:maxinF} \( u \in \F, \) and
    since \( \nu(u,0) = 0 \) Corollary~\ref{cor:2} implies that \(
    \matwo{u}(\{0\}) = 0 .\) 
    
    This also follows from Proposition~\ref{prop:logpot} since \( u
    \geq \lp_{\mu} \) on a neighbourhood of the origin, where \( \mu =
    \sum a_j \delta_{(b_j,0)}.  \)
\end{example}

For now, we seem to be stuck with a rather annoying condition along 
the \( z_{2} \)-axis. When we take away the radially along one of the 
axis we have to restrict the Lelong number along that axis. However 
with little effort we can change that condition to a somewhat weaker, 
or at least more sensible, condition.


\begin{lemma}\label{lem:radial}  
    Assume that \( u \in \F{(D^2)} \) is radial in the first variable,
    \ie\ \( u(|z_1|,z_2) = u(z_1,z_2), \) and let \( \pmb{e_2} =
    (0,1).  \) If there is a constant \( c, \) \( 0 \leq c <
    + \infty \) such that \( \nu(u_{\pmb{e_2}},0) = c, \) and \(
    u(|z_{1}|,z_2) = u(z_1,z_2) \) then 
    \[ 
      \matwo{u}(\{0\}) =
      \big(\ddc{(\max\{c\log|z_2|,u(z_1,z_2)\})}\big)^2 (\{0\}).
    \]
\end{lemma}

\begin{proof}
Set \( h_k = \max(-1,1/k\log|z_1|,1/k\log|z_2|), \) then \( h_k \in
\E_0(D^2), \) and as in the second part of the proof of
Proposition~\ref{prop:cone} we get
\[
   \int_{D^2} -h \, \matwo{u} \geq \int_{D^2}
   -h \, (\ddc{(\max\{c\log|z_2|,u(z_1,z_2)\})}\big)^2,
\]
letting \( k \to \infty \) we get the inequality
\[
   \matwo{u}(\{0\}) \geq
   \big(\ddc{(\max\{c\log|z_2|,u(z_1,z_2)\})}\big)^2
    (\{0\}).
\]
We will proceed to prove an inequality in the opposite direction.

Let us introduce the auxiliary function \( \varphi(\zeta) :=
u(0,\zeta).  \) Note that \( \varphi \in \sh^-(D).  \) Fix \( \zeta
\in D, \) then we have \( \varphi(\zeta) = u(0,\zeta) \leq
\sup_{|z|=r}u(z,\zeta) = u(z,\zeta), \) \( \forall z \in D, \) hence
\begin{equation}\label{eq:trick1}
    u = \max\{\varphi, u\}.
\end{equation}

Using  Riesz
decomposition theorem we can write \( \varphi = p_{\mu} + h, \) where 
\( p_{\mu} \) is the potential of the measure \(  \Delta \varphi \)
and \( h \) is harmonic.

Since \( \nu(u_{\pmb{e_2}},0) = (\Delta \varphi)(\{0\}) = c, \) we
have \( \varphi(\zeta) = c \log|\zeta| + p_{\mu'}(\zeta) + h(\zeta),
\) where the measure \( \mu' \) has no atom at the origin.  By
applying Proposition~\ref{prop:logpot} to the potential \( p_{\mu'} \)
we see that there exists an increasing  sequence of subharmonic functions \( f_{j},
\) with \( f_{j}\nearrow 0 \) (except at the origin), and a decreasing
sequence of numbers \( r_j \searrow 0, \) such that \(  D(r_{j})
\subset \{ \varphi(\zeta) > c\log|\zeta| +f_{j}(\zeta) \} , \) with \( r_{j}
> 0.  \) 

Thus on \( D(1) \times D(r_j) \) we have
\begin{equation}\label{eq:trick2}
    c\log|z_2|+f_j < u(0,z_2) \leq u(z_1,z_2).
\end{equation}

Define a function \[ v_{j}(z_1,z_2) := \max(u(z_1,z_2), c\log|z_2| +f_{j}(z_2) ).
\] Then \( v_j = u \) on \( D(1) \times D(r_{j}),  \)  by the
Equation~\eqref{eq:trick1} and the Inequality~\eqref{eq:trick2}.

Take the sequence \( \{h_k\} \subset \E_0(D^29 \) as above
\[ 
   \int_{D(1) \times D(r_{j})} -h_k \, \matwo{u} = \int_{D(1) \times
   D(r_{j})} -h_k \, \matwo{v_{j}} \leq \int_{D^2} -h_k \, \matwo{v_{j}}.
\]
Since \( v_j \nearrow \max(u(z_1,z_2), c\log|z_2|) \) pointwise
(except at the origin) we 
get according to Theorem~\ref{Thm:CX} that \( \matwo{v_j} \weakto 
\matwo{\max\{u(z_1,z_2), c\log|z_2|\}}. \) In particular
\[ 
   \int_{D^2} -h_k \, \matwo{u} \leq \int_{D^2} -h_k \, \matwo{v_{j}}
   \to \int_{D^2} -h_k \, \matwo{\max\{u(z_1,z_2), c\log|z_2|\}}.
\]
Let \( k \to \infty \) to get \( \matwo{u}(\{0\}) \leq
\matwo{(\max\{c\log|z_2|,u(z,z_2)\})}(\{0\}), \) and the desired
inequality is proved.
\end{proof}

Now we are in position to give the proof of Theorem~\ref{thm:radial2},
which is a generalization of Proposition~\ref{prop:poly-radial}.

\medskip
\noindent\textbf{Proof of Theorem~\ref{thm:radial2}}
Since \( u \in \psh \cap L^{\infty}(D^2\setminus K) \) we have
that \( u \in \F(D^2).  \) Also \( \nu(u_{\pmb{e_2}},0) = c<+\infty,
\) since if \( \nu(u_{\pmb{e_2}},0) = + \infty \) then \( u \equiv
-\infty \) on a neighbourhood of the origin along the \( z_2 \)-axis,
but then \( u \not \in L^{\infty}(\Omega \setminus K).  \) In the same way there
is a \( M>0 \) such that \( u(z,w) > M \log|z|.  \)

According to Lemma~\ref{lem:radial} we have
\[ 
	\matwo{u}(\{0\}) \leq \matwo{(\max\{c\log|w|,u(z,w)\})}(\{0\}).
\]  
Since 
\[
	\nu(\max\{c\log|w|,u(z,w)\},0) \leq \nu(u(z,w),0) = 0, 
\] 
and \(
\max\{M\log|z|,c\log|w|\} \leq \max\{c\log|w|,u(z,w)\}, \)
Proposition~\ref{prop:logestimate} gives
\[
	\matwo{(\max\{c\log|w|,u(z,w)\}(\{0\})} = 0,
\]
thus \( \matwo{u}(\{0\})=0. \) 
\QED

\begin{remark}
    Theorem~\ref{thm:radial2} seems to be a very curious theorem,
    indeed. It may be generalized to functions of the type \(
    u(z_1,\ldots,z_{n-1},z_n) =  u(|z_1|,\ldots,|z_{n-1}|,z_n),  \) by
    a similar reasoning and using the ideas of \cite{Wiklund:04}. However 
    the truly interesting generalization to functions of type
    \( u(z_1,z_2,\ldots,z_n) =  u(z_1,z_2,\ldots,|z_n|), \) which of
    course would answer the main question, is elusive to the author.
\end{remark}

\end{document}